\documentclass[journal]{IEEEtran}
\IEEEoverridecommandlockouts
\pdfoutput=1

\usepackage{calc}

\usepackage{amsmath}
\usepackage{amssymb}
\usepackage{amsthm}
\usepackage{multirow}		
\usepackage{subfig}
\usepackage{stfloats}
\usepackage{xfrac}			
\usepackage{booktabs}
\usepackage{color}

\usepackage{lipsum}
\usepackage{mathtools}
\usepackage{cuted}

\usepackage{url}	

\usepackage{accents}

\allowdisplaybreaks

\newtheorem{Thm}{Theorem}

\usepackage{wrapfig}

\newtheorem{rem}{Remark}
\newtheorem{assumption}{Assumption}

\newcommand{\beq}{\begin{equation}}
\newcommand{\eeq}{\end{equation}}
\newcounter{l1}
\newcounter{l2}
\newcounter{l3}
\newcommand{\bdotlist}{\begin{list}{$\bullet$}{}}
\newcommand{\bboxlist}{\begin{list}{$\Box$}{}}
\newcommand{\bbboxlist}{\begin{list}{\raisebox{.005in}{{\tiny $\blacksquare$ \\ }}}{}}
\newcommand{\bdashlist}{\begin{list}{$-$}{} }
\newcommand{\blist}{\begin{list}{}{} }
\newcommand{\barablist}{\begin{list}{\arabic{l1}}{\usecounter{l1}}}
\newcommand{\balphlist}{\begin{list}{(\alph{l2})}{\usecounter{l2}}}
\newcommand{\bAlphlist}{\begin{list}{\Alph{l2}.}{\usecounter{l2}}}
\newcommand{\bdiamlist}{\begin{list}{$\diamond$}{}}
\newcommand{\bromalist}{\begin{list}{(\roman{l3})}{\usecounter{l3}}}

\usepackage[noadjust]{cite}
\usepackage{filecontents}

\let\dagger\undefined
\newcommand{\dagger}{\mathsf{T}}

\begin{document}
\vspace{-15pt}

\title{Requirements for Interdependent Reserve Types}

\author{\vspace{-8pt}Manuel Garcia\IEEEauthorrefmark{1}, Ross Baldick\IEEEauthorrefmark{3}\vspace{-27pt}
\thanks{\vspace{-15pt}}
\thanks{\IEEEauthorrefmark{1} Sandia National Laboratories \ \  \IEEEauthorrefmark{3} University of Texas at Austin}
\thanks{This work was sponsored by the US Department of Energy's Office of Electricity's Advanced Grid Modeling (AGM) program.  Sandia National Laboratories is a multimission laboratory managed and operated by National Technology \& Engineering Solutions of Sandia, LLC, a wholly owned subsidiary of Honeywell International Inc., for the U.S. Department of Energy’s National Nuclear Security Administration under contract DE-NA0003525. This paper describes objective technical results and analysis. Any subjective views or opinions that might be expressed in the paper do not necessarily represent the views of the U.S. Department of Energy or the U.S. Government.}
\thanks{\vspace{-28pt}}}

\maketitle

\begin{abstract}
As renewable energy penetration increases and system inertia levels drop, primary frequency control is becoming a critical concern in relatively small interconnections such as the Electric Reliability Council of Texas (ERCOT).  To address this problem ERCOT is implementing a number of market rule changes including the introduction of a new Fast Frequency Response (FFR) reserve type to the electricity market.  This FFR reserve type aims to help the traditional Primary Frequency Response (PFR) reserve type in arresting frequency decline in the event of a large generator outage. This paper derives reserve requirements to ensure sufficient reserve to arrest frequency decline before reaching the critical frequency threshold while coupling PFR reserve, FFR reserve, and system inertia.  The general reserve requirement places limits on the amount of PFR reserve that can be provided by each unit based on its ramping capabilities.  Two such limits are derived from first principles and another is proposed that is capable of accommodating the equivalency ratio introduced in previous work.  These PFR reserve limits also provide first principles insight into equivalency ratios, which have only been studied empirically in the past.  High-level insights are provided on a large Texas test case. 

\vspace{-0pt}
\end{abstract}

\vspace{-15pt}
\section{Introduction}
\vspace{-5pt}
The electric power system is experiencing unprecedented penetration levels of wind and solar generation. These inverter-based technologies traditionally do not provide inertia or frequency control services, challenging frequency control~\cite{ratnam2020future,du2017forecast}.  In response, Independent System Operators (ISOs) throughout the United States have introduced new services that aim to improve the frequency response of the system~\cite{xu2016comparison}. The Electric Reliability Council of Texas (ERCOT), the ISO in Texas, has recently reached instantaneous wind penetration levels of 57\% of demand~\cite{ERCOTFactSheetFeb2020}.  ERCOT has introduced a new frequency response service intended to improve primary frequency response, redefined the frequency response reserve products considered in the electricity market, and proposed the introduction of real-time co-optimization~\cite{NPRR581,NPRR863}.  In the context of these new ERCOT rules, this paper derives reserve requirements that ensure sufficient primary frequency response reserve to accommodate a pre-defined generator outage. Though motivated by ERCOT, the models used in this paper are general and are applicable outside of Texas.

Our work focuses on reserve types providing primary frequency control, which intend to arrest frequency decline in the event of a sudden loss of generation and is of critical concern for low inertia systems~\cite{ratnam2020future}.  The derived reserve requirements are unaffected by all other reserve types in ERCOT, which are slower acting and thus do not contribute to arresting frequency decline.  As is common, ERCOT defines its reserve requirements to accommodate the largest possible loss of generation~\cite{pddotnuschel2018frequency}. In the context of primary frequency control, sufficient reserve must be procured to restore power balance before the frequency falls below some \mbox{critical frequency threshold. }

Traditional droop control is provided by synchronous generators and operates within \emph{Primary Frequency Response (PFR)} reserve.  Current practices typically enforce PFR reserve requirements under the assumption that generators are capable of matching the droop reference signal, which is proportional to the frequency deviation and typically incorporates a dead-band. (See Section~\ref{Sec:PFRDef}). However, recent studies show that this assumption does not hold during the transient response of a very large contingency due to turbine governor ramping limitations~\cite{Nesbit2020}.  We address this by differentiating \emph{available PFR reserve} from \emph{nominal PFR reserve}, which intuitively represent the amount of PFR reserve that can be delivered before the critical frequency threshold is met, respectively, with and without considering turbine governor ramping limitations.
 
ERCOT introduced a new reserve type termed \emph{Fast Frequency Response (FFR)} reserve, which enables fast acting devices to act quickly during primary frequency control~\cite{NPRR581}.  Participants include fast acting battery storage and demand curtailment capable of responding within a few voltage frequency cycles.  In contrast to PFR reserve, this new product does not exhibit ramping limitations and is expected to fully deploy nearly instantly when the frequency falls below a specified value, e.g. $59.85$Hz in ERCOT. Throughout this paper we will use this step response definition of FFR reserve, which matches that used by ERCOT but may differ from other definitions~\cite{NERCWhiteFFR}.  Since the response of FFR reserve is nearly instantaneous as defined by ERCOT, these devices are assumed to deliver all reserve before the critical frequency threshold is met. 

ERCOT's real-time co-optimization formulation must include algebraic requirements that couple PFR and FFR reserve. References~\cite{chavez2014governor, trovato2018unit, badesa2019optimal, sokoler2015contingency} derive PFR reserve requirements from first principles without accounting for FFR reserve as defined by ERCOT.  Although \cite{trovato2018unit} and~\cite{badesa2019optimal} use a model that incorporates fixed time delays that vary among generators, their model cannot accommodate FFR reserve, which is deployed at a fixed frequency. Reference \cite{sokoler2015contingency} uses a general model of a generator that is capable of accommodating FFR reserve; however, they utilize a pre-determined frequency trajectory, which would effectively also fix the time that FFR reserve is deployed.  In contrast, we accurately model the FFR reserve as being deployed at a time that varies with the frequency trajectory due to the frequency threshold for activation.

Our main results extend our previous work in~\cite{garcia2019real} and \cite{garcia2019PhD}, which provide a first principles derivation of a requirement that couples PFR and FFR reserve. The proposed reserve requirement has two components.  First, the FFR reserve and available PFR reserve must be sufficient to cover the largest possible loss of generation. Second, limits are placed on the amount of available PFR reserve that can be provided by a synchronous generator based on various parameters including system inertia, size of the contingency being considered, and generators' turbine governor ramp rate (or \emph{governor ramp rate} for brevity).  This paper will derive two such limits from first principles. The first is termed the \emph{rate-based PFR reserve limit} and assumes that each generator's droop response exhibits a constant governor ramp rate that is fixed, similar to the ramping model used in~\cite{chavez2014governor}.  We claim that this model accurately depicts the ramping limits observed in~\cite{Nesbit2020}, which are imposed by thermal and stability limitations of a generator.  The second PFR reserve limit is termed the \emph{proportional PFR reserve limit} and assumes that the governor ramp rates vary proportionally with the amount of nominal PFR reserve allocated to a generator, similar to the ramping model used in~\cite{trovato2018unit} and~\cite{badesa2019optimal}.  This model represents dynamic generator models that mimic a low pass filter on the droop \mbox{reference signal.} 

Reserve requirements coupling PFR and FFR reserve have also been proposed in the context of equivalency ratios, which represent the relative effectiveness of FFR reserve to PFR reserve~\cite{li2018design,liu2018participation}.  This requirement was derived empirically through simulation and simply enforces that the weighted sum of PFR and FFR reserve be larger than an empirically found requirement quantity. By redefining the equivalency ratio as the relative effectiveness of PFR reserve to FFR reserve, we observe that the requirement can be reformulated such that the requirement quantity closely matches the size of the largest possible loss of generation according to empirical data from~\cite{liu2018participation}. Approximating the requirement quantity in this way allows for the equivalency ratio requirement constraint to be enforced using our general framework along with a third type of PFR reserve limit, termed the \emph{equivalency ratio PFR reserve limit}, that is similar in nature to the \emph{proportional PFR reserve limit}. In fact, using these two PFR reserve limits along with various simplifying assumptions, this paper provides a novel first principles derivation of the equivalency ratio, which has only been studied empirically in the literature to date. This result provides first principles insight suggesting that equivalency ratios may vary significantly with total procured FFR reserve, particularly at low inertia levels.

The paper is organized as follows.  Section \ref{Sec:NonLin} presents a general reserve requirement that encompasses all derived reserve requirements.  Section \ref{sec:3Contributors} provides a model of the three main contributors to arresting frequency decline: inertia, PFR reserve, and FFR reserve. This section also explains the offered PFR capacity limits imposed by electricity markets today.  Section \ref{sec:RampRateModels} presents the constant governor ramp rate model of droop control and the rate-based PFR reserve limit.  Section \ref{Sec:PropRate} derives the proportional PFR reserve limit and the equivalency ratio PFR reserve limit.  Under various assumptions, this section also derives the equivalency ratio from first principles.  Section \ref{sec:RTCoOpt} places the reserve requirements into a real-time co-optimization problem and observes high-level trends. 

\vspace{-18pt}
\section{A General Reserve Requirement} \label{Sec:NonLin}
\vspace{-8pt}

The amount of PFR reserve provided by each generator is represented by vector $r\!\in\!\mathbb{R}^n$, where $n$ is the total number of generators.  The amount of PFR reserve provided by generator $i$ is denoted $r_i$.  The amount of FFR reserve provided by each FFR resource is represented by vector $b\!\in\!\mathbb{R}^{\beta}$, where $\beta$ is the total number of FFR resources. The amount of FFR reserve provided by FFR resource $j$ is denoted $b_j$. The vector of ones is denoted $\bold{1}$ and a superscript \!$\dagger$\! represents vector transpose.

The general reserve requirement ensures the system is capable of accommodating a large loss of generation of size $L$, which may include multiple simultaneous unit outages. In the case of ERCOT, $L$ is chosen to be the combined capacity of the two largest generators, which approximately amounts to $2750$MW~\cite{liu2018participation}.  Intuitively, this requirement will additionally ensure sufficient reserve to accommodate any less severe contingencies.  The general reserve requirement is as follows:
\begin{equation}
\bold{1}^{\dagger}r+\bold{1}^{\dagger}b\geq L\label{ResReqMain}	
\end{equation}
In the context of primary frequency control, accommodating a generator outage of size $L$ requires that the voltage frequency trajectory remain above some critical frequency threshold denoted $\omega_{\text{min}}$, which is $59.4$Hz in ERCOT.  This critical frequency threshold represents the point at which firm load begins to disconnect from the system as an emergency precaution to avoid a system wide blackout (See Remark~\ref{buffer1hz}).  This requires power balance to be met before the frequency falls below the $\omega_{\text{min}}$ (See Section~\ref{Sec:SwingEq}). Unlike FFR reserve, governor ramping limitations may limit the amount of PFR reserve that can be deployed before the frequency nadir. With this in mind, the PFR reserve amount appearing in the requirement (\ref{ResReqMain}), denoted $r$, must be available to be deployed before the critical frequency threshold is met and thus we will hence-forth refer to this quantity as the \emph{available PFR reserve}.  In contrast, the \emph{nominal PFR reserve}, represented by vector $R\!\in\!\mathbb{R}^n$, will refer to the amount of headroom each generator is maintaining for the purpose of PFR reserve, some of which may not be \emph{available} before the critical frequency threshold is met.  The nominal PFR reserve provided by generator $i$ is denoted $R_i$.

The distinction between nominal PFR reserve and available PFR reserve is essential to understanding and analyzing the role and value of PFR and FFR reserve. Nominal PFR reserve is limited by the generator's head-room as well as the generator's offered PFR capacity (See Section \ref{Sec:PFRDef}).  Available PFR reserve is limited by nominal PFR reserve, e.g. $r_i\!\leq\! R_i$, because a generator's droop response should not be expected to exceed the amount of headroom it has procured for the purpose of PFR reserve. Available PFR reserve may also be limited by other factors including the generator's turbine governor ramping ability and the time taken to reach the critical frequency threshold.  The remainder of this paper is dedicated to deriving different limits that can be placed on the available PFR reserve to ensure it can be delivered before the critical frequency threshold is met. 

\vspace{-5pt}
\rem{\normalfont Reserve requirements in ERCOT are designed to maintain the frequency above $\omega_{\text{min}}=59.4$Hz and firm load begins to disconnect at $59.3$Hz~\cite{liu2018participation}. This $0.1$Hz margin accommodates potential errors, for example, error in frequency measurement that may occur during transient conditions.\label{buffer1hz}}\normalfont

\vspace{-14pt}
\section{Three Contributors to Arresting Frequency}\label{sec:3Contributors}
\vspace{-5pt}
The three main contributors to arresting frequency decline in response to a large generator outage are inertia, PFR reserve, and FFR reserve.  Each of these will now be modeled in detail.  

\vspace{-5pt}
\subsection{Inertia and Frequency Dynamics}\label{Sec:SwingEq}
\vspace{-2pt}
Voltages in the system are modeled as quasi-steady state sinusoids whose frequency may be slowly varying.  Moreover, this voltage frequency at time $t$ is modeled as being the same at each generator in the system and is denoted $\omega(t)$.  The total post-outage inertia is $M$ (in units of Watt-seconds or Ws) and represents the sum of inertia values for all generators still in service after the outage. The system dynamics are represented by the \emph{swing equation}~\cite{anderson2008power}, which is expressed as follows: 
\begin{equation}
{\textstyle \frac{d\omega(t)}{dt}}={\textstyle \frac{\omega_0}{2M}}(\bold{1}^{\dagger}m(t)-e(t)),\label{SwingEqn}	
\end{equation}
where $m(t)\in\mathbb{R}^n$ represents the vector of mechanical power output from the turbine governor of each generator in the system and $e(t)\in\mathbb{R}$ represents the total net electrical demand in the system.  This model makes the simplifying assumption that there is no system damping.  The \emph{nominal frequency} is denoted $\omega_0$ and will be assumed to be the frequency just prior to the time of the generator outage.

\vspace{-5pt}
\subsection{Primary Frequency Response Reserve and Droop Control}\label{Sec:PFRDef}
PFR reserve is intended to be compatible with conventional generator droop control, which increases the mechanical power output of the turbine governor $m_i(t)$ in response to a large generator outage.   We conservatively assume that generators contracting to provide PFR reserve are the only generators that actually do provide droop control.  This assumption deviates slightly from ERCOT requirements, which instead widen the droop control dead-band for all generators not contracting to provide PFR reserve~\cite{NPRR863}.  PFR reserve is provided by generators that may also be selling power into the electricity market.  The nominal PFR reserve $R_i$ must satisfy \mbox{$G_i+R_i\leq \bar{G}_i$,} where $G_i$ is the dispatched electric power generation of generator $i$ and $\bar{G}_i$ is its capacity. Furthermore, each generator has an \emph{offered PFR capacity} denoted $\bar{R}_i$.  With this in mind, the private constraints for all generators are written as follows: 
\begin{equation} \underline{G}_i\leq G_i \leq\bar{G}_i-R_i \text{ and }  0\leq R_i\leq \bar{R}_i \ \ \forall i\in[1,\ldots,n].\label{PrivateGenerators}	
\end{equation}

The ISO typically has qualification requirements that enforce a limit on the offered PFR capacity $\bar{R}_i$ that a generator can offer into the market.  To derive this limit on the offered PFR capacity $\bar{R}_i$ we first must describe standard droop control in detail. Subsequently, the concept of available PFR reserve will be detailed, which further limits the amount of PFR that can be deployed in particular system conditions, to account for turbine governor ramping limitations.

\subsubsection{Standard Droop Control}\label{StandardDroop}
\begin{figure}[h]
\vspace{-8pt}
\centering
\includegraphics[scale=.5]{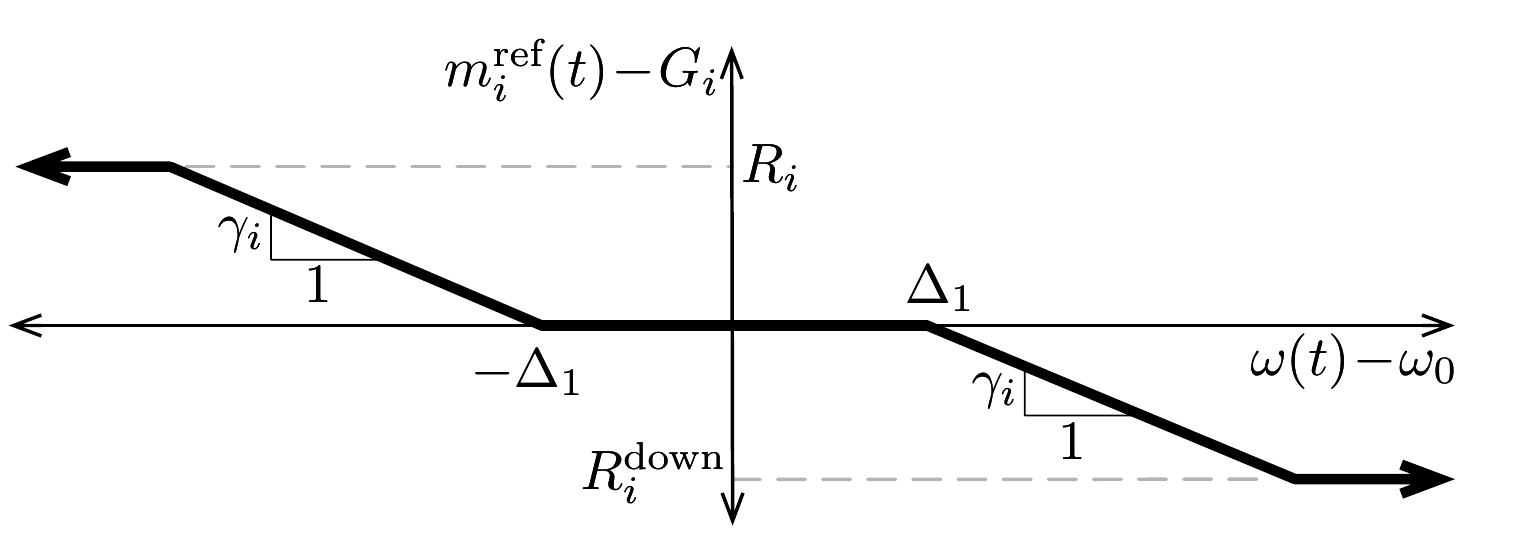}
\vspace{-10pt}
\caption[Droop reference signal with dead-band.]{\label{Fig:DroopSig}Droop reference signal with dead-band.}
\vspace{-5pt}
\end{figure}

Generators providing PFR reserve respond to local frequency via droop control by adjusting the reference mechanical power output of their turbine governor $m_i^{\text{ref}}(t)$ on the time scale of $\frac{1}{\omega_0}\approx 0.016$ seconds.  In the context of droop control, the generation value $G_i$ represents the nominal value of $m_i^{\text{ref}}(t)$ around which the adjustments are made and is updated each time the real-time market clears, which occurs every five minutes. Furthermore, during droop control the reference mechanical power output of each generator's turbine governor is limited implicitly by the need to preserve capacity for reserves, is adjusted depending on the generator's local frequency deviation, and has a dead-band of $\Delta_1:=\omega_0-\omega_1$ where \mbox{$\omega_1<\omega_0$} represents the low end of the dead-band.  This reference signal, or \emph{droop reference signal}, is illustrated in Figure~\ref{Fig:DroopSig} where $\gamma_i$ is the droop constant for generator $i$ and $R^{\text{down}}_i$ represents the nominal down PFR reserve, which is not detailed in this work but is analogous to $R_i$.  Notice that the droop reference signal is limited by the nominal PFR reserve $R_i$.  This is because the nominal PFR reserve $R_i$ typically matches either the PFR reserve capacity $\bar{R}_i$ or the headroom of the generator $\bar{G}_i-G_i$.  If this limit is not imposed, then the droop reference signal may incorrectly instruct the generator to produce more than its capacity $\bar{G}_i$.

The reference mechanical power output then traverses the turbine governor dynamics of the generator to produce the realized mechanical power output, denoted $m_i(t)$.  These turbine governor dynamics can be very complicated and are not detailed in our work. (See~\cite{peydayesh2017simplified} and~\cite{kundur1994power} for details).  

\subsubsection{Offered PFR Capacity Limits}\label{PFROfferCaps}

The droop reference signal $m_i^{\text{ref}}(t)$ cannot be larger than $G_i+\gamma_i(\omega_0-\omega_{\min})-\gamma_i\Delta_1$ without the critical frequency threshold $\omega_{\min}$ being violated.  For this reason, ISOs should impose the following offered PFR capacity limit for each generator $i$:
\begin{equation}
\bar{R}_i \leq \gamma_i(\omega_0-\omega_{\min})-\gamma_i\Delta_1\label{Constraint:TradLim}
\end{equation}
The droop constant $\gamma_i$ is chosen based on a required droop percentage imposed by the ISO.  The droop percentage represents the percent change in frequency required to achieve a governor change of 100 percent capacity. Let $\nu_i$ represent the droop percentage for each generator expressed as a fraction.  For example, in ERCOT the droop percentage is typically set to $5\%$ and so $\nu_i=0.05$. The proportionality droop constant $\gamma_i$ used during droop control satisfies \mbox{$\gamma_i\nu_i\omega_0-\gamma_i\Delta_1=\bar{G}_i$} and can be determined as follows:
\begin{equation}\gamma_i=\tfrac{\bar{G}_i}{\nu_i\omega_0-\Delta_1}\label{DroopConstantEqn}\end{equation}
This definition of the droop constant is consistent with BAL-001-TRE-1, the reliability standard that details primary frequency response requirements in ERCOT~\cite{Bal001}. Following from (\ref{Constraint:TradLim}) and (\ref{DroopConstantEqn}) the offered PFR capacity limit is written as follows:
\begin{equation}
\bar{R}_i \leq \tfrac{\bar{G}_i(\omega_0-\omega_{\min}-\Delta_1)}{\nu_i\omega_0-\Delta_1}\approx \tfrac{\bar{G}_i(\omega_0-\omega_{\min})}{\nu_i\omega_0}.\label{ERCOTLimApprox}
\end{equation}
The approximation assumes that the dead-band for droop control $\Delta_1$ is very small.  Although the approximation over estimates the offered PFR capacity limit, the approximation error is typically very small and is easily accommodated by the conservatively chosen critical frequency threshold $\omega_{\text{min}}$ (See Remark~\ref{buffer1hz}).  In fact, the typical dead-band in ERCOT is $\Delta_1=0.017$Hz, which is significantly smaller than the value $\omega_0-\omega_{\min}=0.6$Hz and the typical value of $\nu_i\omega_0=3$Hz for the typical droop percentage of $5\%$.  Furthermore, the critical frequency threshold is $\omega_{\text{min}}=59.4$Hz in ERCOT and ERCOT uses the approximation outlined in (\ref{ERCOTLimApprox}).  In this case the offered PFR capacity limit is $0.2\bar{G}_i$ for a generator $i$ with $5\%$ droop.  This is consistent with ERCOT protocols~\cite{NPRR863}.

\subsubsection{Available PFR Reserve and Ramping Limitations}

The reference mechanical power output must traverse the turbine governor dynamics of the generator to produce realized mechanical power output, denoted $m_i(t)$.  A salient feature of these turbine governor dynamics is that the mechanical power output $m_i(t)$ tends to lag the reference mechanical power input $m^{\text{ref}}_i(t)$, particularly if the reference signal changes quickly.  For this reason, it is possible that a generator's PFR reserve is not fully available before the critical frequency threshold is met, effectively exhibiting ramp limitations that restrict its output.  Current practices do not explicitly account for these ramping limitations.  We address this shortfall in current practices by differentiating between nominal PFR reserve $R_i$ and \emph{available} PFR reserve $r_i$, which represents the amount of PFR reserve that is actually available as increased generation before the critical frequency threshold is reached.  

 \vspace{-7pt}
\subsection{Fast Frequency Response Reserve}\label{Sec:FFRDef}
\vspace{-2pt}
We assume that FFR reserve can be fully deployed instantaneously and can be sustained for several minutes, until slower acting reserve is capable of responding (See Remark~\ref{Rem:InstantaneousFFR}).  FFR reserve can be provided by any device that meets these requirements.  FFR capable devices include fast-acting battery storage and load-shedding.  The amount of FFR reserve provided by FFR resource $j$ is denoted $b_j$ and FFR resource $j$ provides an \emph{offered FFR capacity} of $\bar{b}_j$.  The private constraints for all FFR resources are written as follows:
\begin{equation}
0\leq b_j\leq\bar{b}_j \ \ \ \forall j \in [1,\ldots,\beta].\label{PrivateBatteries}	
\end{equation}

The FFR reserve is deployed when the frequency drops below a frequency threshold of \mbox{$\omega_2<\omega_1$}, where $\omega_1$ is the frequency corresponding to the dead-band of droop control.  Note that $\omega_2$ is typically significantly lower than the frequency $\omega_1$.  In fact, FFR reserve is deployed only during emergencies involving the largest generator outages as opposed to PFR reserve, which is used for essentially all contingencies.  When deployed, the FFR reserve instantaneously decreases the net electrical demand in the system $e(t)$ by an amount $\bold{1}^{\dagger}b$.  We additionally introduce the non-negative constant $\Delta_2\!:=\!\omega_1\!-\omega_2$.

\rem{\normalfont ERCOT requires FFR reserve to be capable of being completely deployed within $0.25$ seconds of being called upon and being sustained for at least $15$ minutes~\cite{NPRR863}.  A simple extension of our FFR model can include this time delay of $0.25$ seconds before the instantaneous response occurs.  This extension can be easily accommodated and would require small adjustments in the proof of Theorem~\ref{Thm:RateBasedPFRLimit}, which are not pursued in this paper due to space constraints.\label{Rem:InstantaneousFFR}}\normalfont

\vspace{-5pt}
\section{Rate-Based PFR Reserve Limits}\label{sec:RampRateModels}
\vspace{-2pt}
This section presents the rate-based PFR reserve limit originally derived in previous work~\cite{garcia2019real,garcia2019PhD}.

\vspace{-8pt}
\subsection{Simple Turbine Governor Model}\label{Sec:TwoModels}
\begin{figure}[h]
\vspace{-8pt}
\centering
\includegraphics[scale=.45]{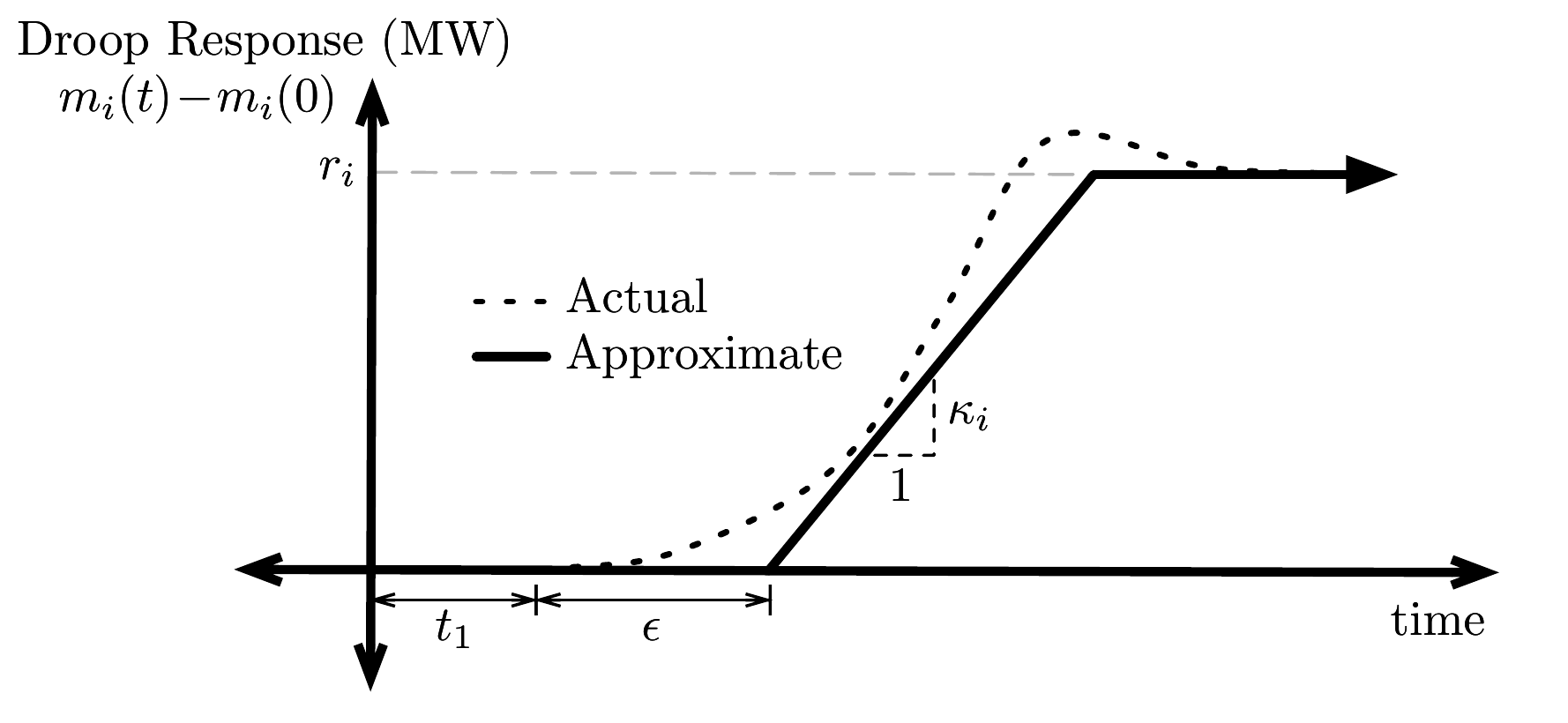}
\vspace{-10pt}
\caption[Turbine governor response to generator outage.]{\label{GovResp} Turbine governor response to generator outage. This figure is based on figure 2 in~\cite{chavez2014governor}.}
\vspace{-10pt}
\end{figure}

We need only characterize each generator's turbine governor response to the very specific situation where a large loss of generation occurs since all smaller contingencies will result in smaller frequency excursions.  Such a response is similar to that of a step response in the reference mechanical power output because of the fast frequency drop.  This type of response is illustrated in \mbox{Figure~\ref{GovResp}}.  The approximate piecewise linear model shown in this figure is adopted from~\cite{chavez2014governor}. 

The mechanical power output of generator $i$ is assumed to match its dispatched generation $G^d_i$ at the time of the generator outage $t\!=\!0$, e.g. $m_i(0)\!=\!G^d_i$.  Following the outage at $t\!=\!0$ the frequency begins to drop.  At the time $t\!=\!t_1$ the frequency reaches the lower end of the frequency dead-band $\omega(t_1)\!=\!\omega_0-\Delta_1$.  Subsequently, the turbine governor is modeled as experiencing a small time delay $\epsilon$.  Although this time delay is assumed the same for all generators, this is only a simplifying assumption that can be easily extended. After exhibiting this time delay, the mechanical power output is modeled as having a constant \emph{governor ramp rate} $\kappa_i$ that continues until all available PFR reserve $r_i$ is deployed.

\vspace{-10pt}

\subsection{Sufficient Condition for Satisfying Frequency Threshold}\label{Sec:SuffCond}

References \cite{garcia2019real,garcia2019PhD} prove that the rate-based PFR reserve limit serves as a sufficient condition for adequate reserve procurement.  In order for this result to hold, we must impose two assumptions regarding the response of the system.  First, we assume there is enough reserve to restore power balance, which is ensured by the general reserve requirement (\ref{ResReqMain}).  Second, we assume the FFR reserve is deployed after the PFR reserve begins ramping upward.  This assumption is reasonable because, as discussed above, the PFR reserve dead-band threshold $\Delta_1$ is much tighter than the FFR reserve dead-band threshold $\Delta_1\!+\!\Delta_2$ and the delay $\epsilon$ is typically small.   

\assumption{\normalfont We assume there is sufficient reserve to restore power balance, which is ensured by the general reserve requirement (\ref{ResReqMain}).  We additionally assume that the FFR reserve is deployed after the PFR reserve begins to ramp upward, which is satisfied if the post-contingency inertia satisfies:
\begin{align}
	&M\geq{\textstyle \frac{\epsilon L \omega_0}{2\Delta_2}}\label{ass02}
\end{align}

\vspace{-0pt}
\label{ass0}}\normalfont
The thresholds are set according to the ERCOT NPRR 863~\cite{NPRR863}, resulting in \mbox{$\omega_0=60$Hz,} \mbox{$\omega_1=59.9833$Hz,} and \mbox{$\omega_2=59.85$Hz.} These parameters will be used in all numerical results in this paper along with a PFR time delay of $\epsilon=0.2$ seconds.   Furthermore, $L$ is typically set to $2750$MW to represent the loss of the two largest generators in ERCOT.  With these parameters, the assumption (\ref{ass02}) can be interpreted as a lower bound on the system inertia of approximately $123$GWs.  

The following theorem presents the rate-based PFR reserve limit, where the sum of all FFR reserve is denoted $\tilde{b}:=\bold{1}^{\dagger}b$ and a constant is introduced as $\Delta_3:=\omega_0-\Delta_1-\Delta_2-\omega_{\text{min}}$.

\Thm{\label{Thm:RateBasedPFRLimit}\hspace{-3pt}Under Assumption~\ref{ass0},  the frequency $\omega(t)$ will remain above the minimum frequency threshold $\omega_{\text{min}}$ if:
\normalsize{
\begin{equation}
r_i\leq \kappa_i h(M,\bold{1}^{\dagger}b) \ \ \ \forall i\in[1,\ldots,n]\label{ResReqRamp}
\end{equation}}
where the limit function $h(M,\tilde{b})$ is expressed as follows:
\begin{equation}h(M,\tilde{b}):=\tfrac{\tfrac{4M}{\omega_0}(\Delta_2+\Delta_3-\frac{\omega_0}{2M}\epsilon L)^2(L-\tilde{b})}{\left(\hspace{-1pt}\tilde{b}\sqrt{\Delta_3}\hspace{-1pt}-\hspace{-1pt}\sqrt{(\Delta_2\hspace{-1pt}+\hspace{-1pt}\Delta_3\hspace{-1pt}-\hspace{-1pt}\frac{\omega_0}{2M}\epsilon L) L^2\hspace{-1pt}-\hspace{-1pt}(\Delta_2\hspace{-1pt}-\hspace{-1pt}\frac{\omega_0}{2M}\epsilon L)\tilde{b}^2}\right)^2}\label{hFuncDef}\end{equation}

}\normalfont

\proofname: Provided in~\cite{garcia2019real} and~\cite{garcia2019PhD}. In addition to \mbox{Assumption \ref{ass0},} \cite{garcia2019real} and~\cite{garcia2019PhD} assume the frequency falls below $\omega_2$ in response to the outage of size $L$ and the FFR reserve deployment does not \emph{overshoot} causing a positive power imbalance. When these additional assumptions are violated the result still holds because the frequency will still not fall below $\omega_2>\omega_0$.
\qed

 The rate-based PFR reserve limit from (\ref{ResReqRamp}) can be interpreted as a condition that the available PFR reserve must satisfy.  With this in mind, any additional nominal PFR reserve in excess of this limit cannot be utilized before the critical frequency threshold is met. Intuitively, any nominal PFR reserve $R_i$ that exceeds the amount $\kappa_i h(M,\bold{1}^{\dagger}b)$ is not considered \emph{available}.

\begin{figure}[h]
\vspace{-12pt}
\centering
\includegraphics[scale=.4]{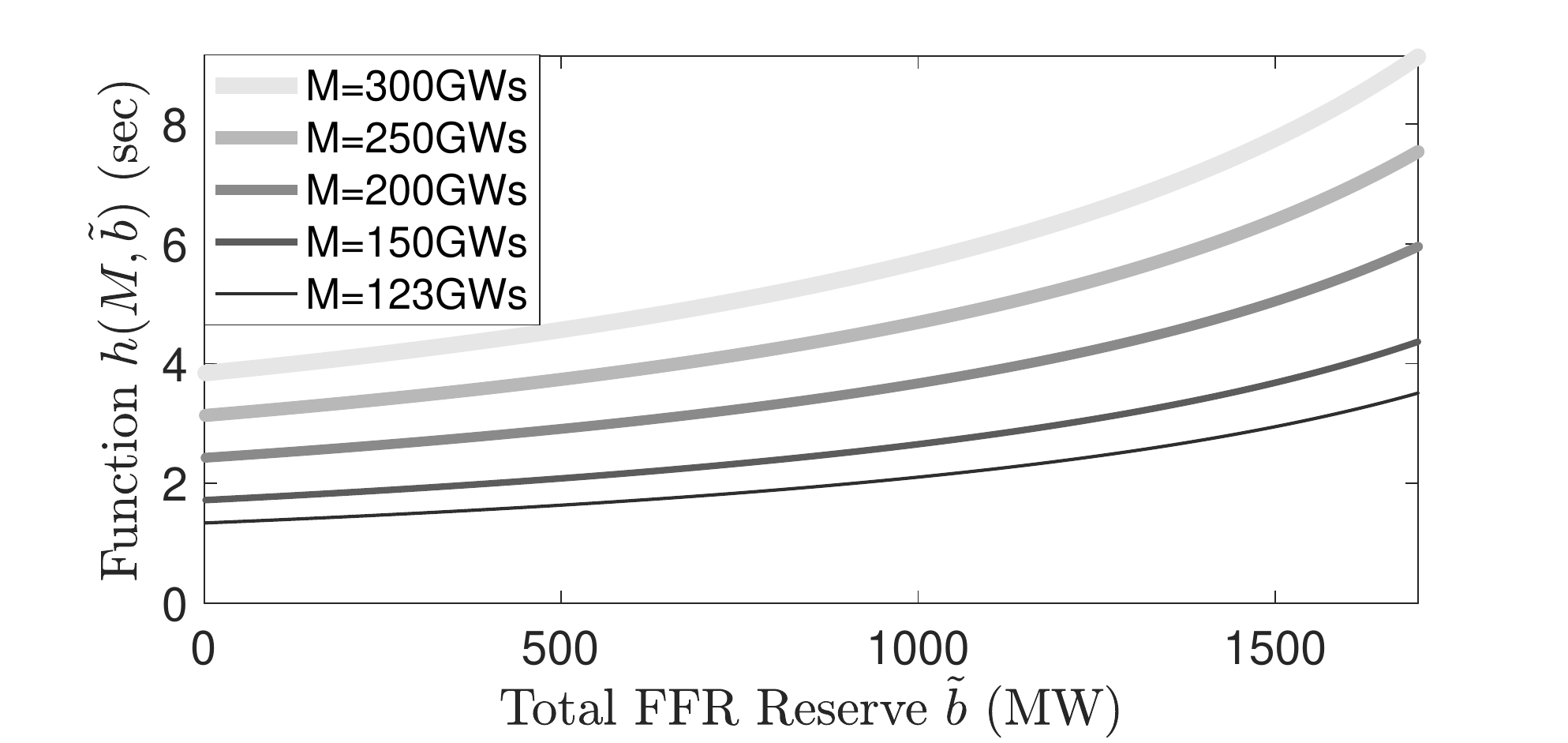}
\vspace{-5pt}
\caption[]{\label{fig:LimitFunc} Function $h(M,\tilde{b})$ with parameters from Section~\ref{Sec:SuffCond}.  Reference~\cite{garcia2019real} originally published this plot with incorrect scaling of the vertical axis. That error is corrected in this plot.}
\vspace{-10pt}
\end{figure}

The function $h(M,\tilde{b})$ is convex in its second argument. Figure~\ref{fig:LimitFunc} provides example plots of $h(M,\tilde{b})$ versus its second argument $\tilde{b}$ for several different values of inertia $M$.  Notice that this function is positive and increasing in $M$ and $\tilde{b}$.  As a result constraint (\ref{ResReqRamp}) allows the nominal PFR reserve for a generator to increase if the system inertia increases, the total FFR reserve increases, or its governor ramp rate $\kappa_i$ increases.

\vspace{-8pt}
\section{Proportional PFR Reserve Limits}\label{Sec:PropRate}
\vspace{-3pt}

This section provides limits to be placed on a generator's available PFR reserve, $r_i$, that increase proportionally with its procured nominal PFR reserve, $R_i$.  Section~\ref{Sec:PropModel} first introduces a turbine governor model that assumes the governor ramp rate $\kappa_i$ is proportional to the nominal PFR reserve $R_i$.  This model results in the \emph{proportional PFR reserve limit}.  Section~\ref{sec:EqRatReq} then provides the equivalency ratio requirement from~\cite{liu2018participation} and explains how this requirement can also be represented using the \emph{equivalency ratio PFR reserve limit}.  Section~\ref{Sec:ApproxEquivRat} then uses these two limits to provide first principles insight into the behavior of equivalency ratios.  

\vspace{-7pt}
\subsection{Proportional Ramp Rate Model}\label{Sec:PropModel}
In the case of large contingencies, the governor ramp rate of a generator, $\kappa_i$, may be considered fixed as suggested in \cite{chavez2014governor}.  Such a model would be accurate if the generator experiences thermal and stability limitations described in~\cite{Nesbit2020}.  

On the other hand, if a generator is not limited in this way, then its governor ramp rate, $\kappa_i$, may increase with its nominal PFR reserve, $R_i$.  This is because the frequency drop in response to a large generator outage is so fast that the droop reference signal as described in Section~\ref{StandardDroop} will be similar to a step change, where the size of the step is equal to the nominal PFR reserve $R_i$.  Assuming a low-pass filter type response from the droop reference signal to the generator mechanical output power, the ramp rate of the turbine governor's mechanical power output will increase with the size of the step reference input to the turbine governor.  

To emulate this dependency on nominal PFR reserve, we model the generator's governor ramp rate as being proportional to the nominal PFR reserve $R_i$, as suggested in~\cite{trovato2018unit} and~\cite{badesa2019optimal}.  Let's denote the proportionality constant by $\lambda_i\geq 0$, which may vary among generators.  Then the governor ramp rate of the generator is $\kappa_i=\lambda_i R_i$. Notice that the nominal PFR reserve amount $R_i$ is changed on 5-15 minute time scales corresponding to the real-time market clearing.  On the other hand, primary frequency control occurs on time scales of tenths of a second.  As a result, the nominal PFR reserve amount and the governor ramp rate $\kappa_i$ can be considered constant on the time scales that primary control is being performed.  With this in mind, the same rate-based PFR reserve limit (\ref{ResReqRamp}) can be enforced as:
\begin{equation}
r_i\leq \lambda_i R_i h(M,\bold{1}^{\dagger}b) \ \ \ \forall i\in[1,\ldots,n].\label{Constraint:ConstTime}
\end{equation}
We will refer to this limit as the \emph{proportional PFR reserve limit}.  Intuitively, the available PFR reserve $r_i$ is limited to be a fraction of the nominal PFR reserve $R_i$, where the fraction $\lambda_i h(M,\bold{1}^{\dagger}b)$ is non-negative and is typically less than one.

\vspace{-5pt}
\subsection{Equivalency Ratio Requirement}\label{sec:EqRatReq}

Reference~\cite{liu2018participation} uses the \emph{equivalency ratio reserve requirement}:
\begin{equation}
\bold{1}^{\dagger}R +\alpha(M)\bold{1}^{\dagger}b\geq \upsilon(M)\label{PengweDuReq}
\end{equation}
where $\alpha(M)$ is termed the \emph{equivalency ratio} and $\upsilon(M)$ is the \emph{frequency response reserve requirement (Rfrr)}, both of which are functions of the total system inertia.  Reference~\cite{liu2018participation} determines these two functions empirically based on simulation studies so that satisfaction of constraint (\ref{PengweDuReq}) will ensure sufficient reserve to prevent the frequency from violating the minimum frequency threshold of $\omega_{\text{min}}=59.4$Hz in response to an outage of size $L$, where $L=2750$MW in their work.  The first three columns of Table~\ref{DuTable} replicate the data from~\cite[table I]{liu2018participation} and show the values of $\alpha(M)$ and $\upsilon(M)$ for different values of inertia $M$.

\begingroup
\setlength{\tabcolsep}{1pt} 
\renewcommand{\arraystretch}{1} 
\vspace{-0pt}
\begin{table}[h!]
\caption[]{\label{DuTable} Parameters appearing in the equivalency ratio requirement from~\cite{liu2018participation}.  The equivalency ratio $\alpha(M)$ and reserve requirement $\upsilon(M)$ are shown for different inertia levels. }
\vspace{-10pt}
\label{tt}
\begin{center}
\begin{tabular}{|c|c|c|c|} 
	\hline
	Total Inertia $M$& Rfrr $\upsilon(M)$&Equivalency &Ratio $\frac{\upsilon(M)}{\alpha(M)}$\\
	(GWs)&(MW)&Ratio $\alpha(M)$&(MW)\\
    \hline
    \hline
    120&5200&2.2&2363.6\\
    \hline
    136&4700&2.0&2350.0\\
    \hline
    152&3750&1.5&2500.0\\
    \hline
    177&3370&1.4&2407.1\\
    \hline
    202&3100&1.3&2384.6\\
    \hline
    230&3040&1.25&2432.0\\
    \hline
    256&2640&1.13&2336.3\\
    \hline
    278&2640&1.08&2444.4\\
    \hline
    297&2240&1&2240.0\\
    \hline
\end{tabular}
\end{center}
\vspace{-21pt}
\end{table}

\endgroup

Table~\ref{DuTable} additionally presents the ratio $\frac{\alpha(M)}{\upsilon(M)}$, which is approximately constant across all inertia levels.  We suggest that the small 10\% variation in these ratios is likely due to process noise associated with the empirical simulation process.  Furthermore, this ratio is approximately equal to but slightly less than the magnitude of the outage being accommodated $L\!=\!2750$.  This slight mismatch may be due to contributions of frequency responsive load in the simulations of reference~\cite{liu2018participation} (See Remark~\ref{Rem:FreqRespLoads}).  With this caveat, we suggest that this observation merits the approximation of $\upsilon(M) \approx \alpha(M)L$.  For this reason, the constraint (\ref{PengweDuReq}) can be approximated as the following \emph{reformulated equivalency ratio requirement}:
\vspace{-3pt}
\begin{equation}
\tfrac{1}{\alpha(M)}\bold{1}^{\dagger}R +\bold{1}^{\dagger}b\geq L\label{PengweDuReq2}
\end{equation}
\vspace{-17pt}

The \emph{reformulated equivalency ratio}, denoted $\tfrac{1}{\alpha(M)}$, represents the relative effectiveness of PFR to FFR reserve. Furthermore, in the context of our general reserve requirement (\ref{ResReqMain}), this reformulation suggests that the following \emph{equivalency ratio PFR reserve limit} be placed on the available PFR reserve:
\begin{equation}
r_i\leq\tfrac{1}{\alpha(M)}R_i \ \ \ \forall i\in[1,\ldots,n].\label{EquivRatioPFRLimit}
\end{equation}
In fact, the reformulated equivalency ratio reserve requirement (\ref{PengweDuReq2}) holds if and only if there exists a vector $r\in\mathbb{R}^n$ such that constraints (\ref{EquivRatioPFRLimit}) and (\ref{ResReqMain}) hold.

\rem{\normalfont \hspace{-2pt}Some load types (e.g. induction motors) will naturally contribute to primary frequency control without explicitly providing reserve~\cite{ERCOTManual,ERCOTWhiteInertia}. As a result, setting the reserve requirement $L$ equal to the contingency size (e.g. $2750$ MW in ERCOT) is a conservative approximation.  In the context of the reformulated equivalency ratio requirement, the requirement quantity $L$ is represented by the ratio column in Table \ref{DuTable}.  This column shows that the requirement quantity $L$ can be adjusted downward from $2750$MW to approximately $2500$MW.\label{Rem:FreqRespLoads}}\normalfont

\vspace{-10pt}
\subsection{Proportional Ramp Rate Model and Equivalency Ratios}\label{Sec:ApproxEquivRat}
\vspace{-2pt}

Although~\cite{liu2018participation} initially justified the use of equivalency ratios empirically, our analysis provides insight into equivalency ratios established from first principles.  Specifically, the equivalency ratio PFR reserve limit from (\ref{EquivRatioPFRLimit}) and the proportional PFR reserve limit from (\ref{Constraint:ConstTime}) are both proportional to the nominal PFR reserve $R_i$, suggesting that the proportionality constants may be similar.  With this in mind, under the assumption that all proportionality constants $\lambda_i$ are approximately the same, the equivalency ratios can be approximated as follows:
\vspace{-5pt}
\begin{equation}
\alpha(M)\approx\tfrac{1}{\lambda_i h(M,\bold{1}^{\dagger}b)}\label{EquivRatioApprox}
\end{equation}
\vspace{-15pt}

\begin{figure}[h]
\centering
\vspace{-10pt}
\includegraphics[scale=.4]{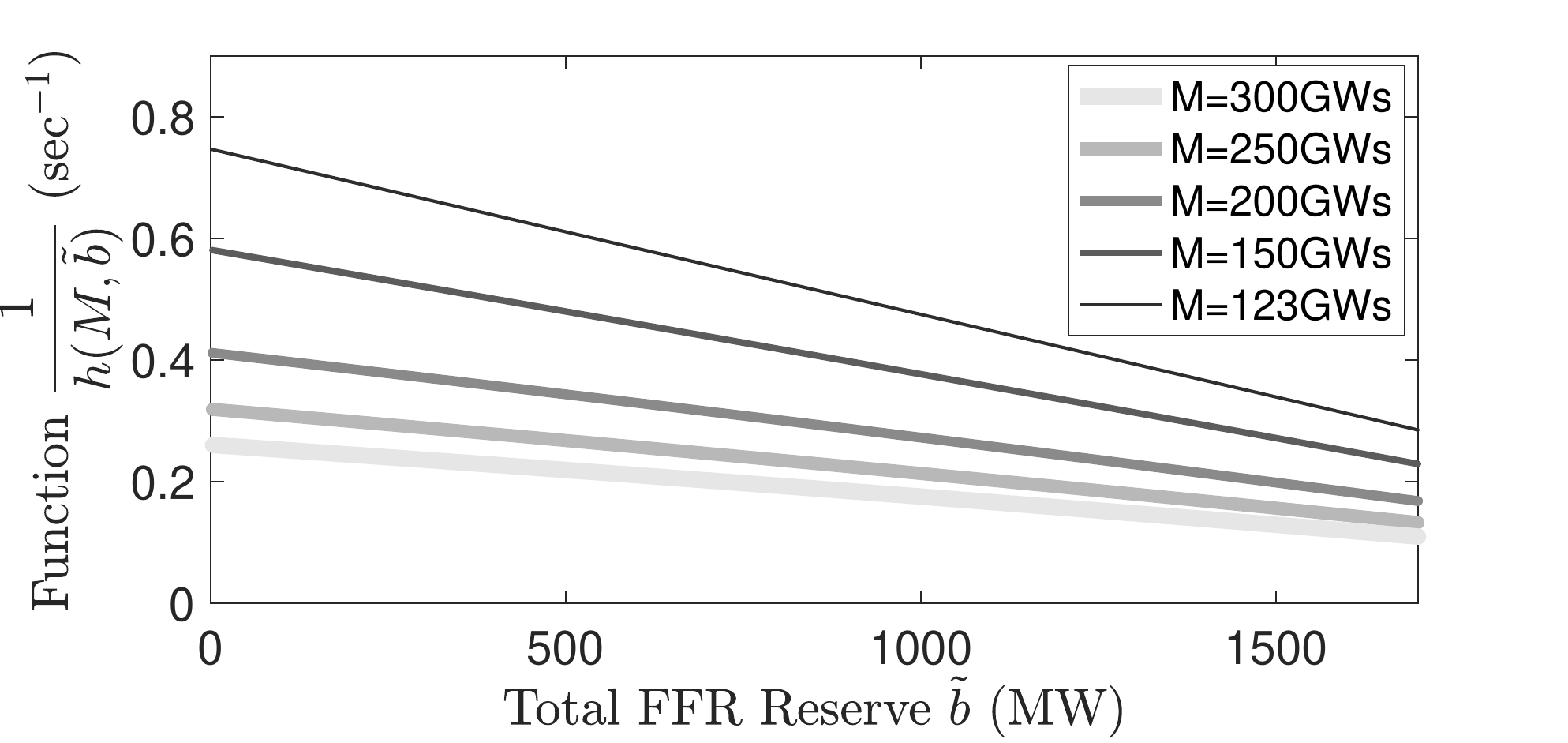}
\vspace{-5pt}
\caption[Inverse of the limit function as a function of total FFR reserve.]{\label{fig:LimitFuncInv} Function $\tfrac{1}{h(M,\bold{1}^{\dagger}b)}$ with parameters from Section~\ref{Sec:SuffCond}.}
\vspace{-5pt}
\end{figure}

To better understand this approximation, Figure~\ref{fig:LimitFuncInv} plots the function $\tfrac{1}{h(M,\bold{1}^{\dagger}b)}$ versus the total FFR reserve $\bold{1}^{\dagger}b$ for different values of inertia $M$. This function varies only slightly in the total FFR reserve argument at high inertia levels.  In other words, the function $\tfrac{1}{h(M,\bold{1}^{\dagger}b)}$ can be reasonably approximated as being constant in the total FFR reserve $\bold{1}^{\dagger}b$ when the inertia level is high, further justifying the approximation (\ref{EquivRatioApprox}).  

Approximation (\ref{EquivRatioApprox}) provides important insight into equivalency ratios from first principles.  Specifically, equivalency ratios are appropriately approximated as being constant in the total FFR reserve, $\bold{1}^{\dagger}b$, at the high inertia levels that are similar in magnitude to the inertia levels seen today in ERCOT.  However, as the inertia levels drop, the equivalency ratio may vary significantly with the total FFR reserve.  As is seen in Figure~\ref{fig:LimitFuncInv}, the slope of the equivalency ratio with respect to the total FFR reserve is approximately three times larger when the inertia is $123$GWs as opposed to $300$GWs. Future work will extend the empirical results from~\cite{liu2018participation} to observe the dependency of equivalency ratios on the total FFR reserve.

\vspace{-10pt}
\section{Real-Time Co-Optimization}\label{sec:RTCoOpt}
\vspace{-6pt}

This section places reserve requirements into simplified real-time co-optimization problems to identify high-level tendencies of the proposed PFR reserve limits and to compare them to the previously proposed equivalency ratio reserve requirement.  

\vspace{-15pt}
\subsection{Texas Test Case}
\vspace{-6pt}
We study a 2000 bus test case that represents the ERCOT system and is based on publicly available data~\cite{birchfield2017grid,xu2018modeling}.  The cost for generator $i$ is represented by the convex cost function $C_{i}(G_{i})$.  For simplicity, we assume the offer price for PFR and FFR reserve is zero and as a result the optimal total procured FFR reserve $\bold{1}^{\dagger}b$ matches the total offered FFR capacity which will be fixed at $\bold{1}^{\dagger}\bar{b}=600$MW.  The 50 natural gas generators with the largest capacity are selected to provide PFR reserve and their offered PFR capacity is set to $\bar{R}_i=0.2 \bar{G}_i$, consistent with Section~\ref{PFROfferCaps}. We consider a loss of generation in the amount of $2750$MW, which represents the two largest nuclear plants in Texas. To adjust for frequency responsive loads, the parameter $L$ is set to $2500$MW as explained in Remark \ref{Rem:FreqRespLoads}. All other parameters are consistent with previous sections.

\begingroup
\setlength{\tabcolsep}{4pt} 
\renewcommand{\arraystretch}{1.25} 
\vspace{-10pt}
\begin{table*}[t!]
\vspace{-10pt}
\caption[]{\label{Table:OptProbs} Three co-optimization problems differentiated by their reserve requirements only.}
\vspace{-10pt}
\label{tt}
\begin{center}

\begin{tabular}{|c|c|c|}
\hline
Equivalency Ratio Requirement&Rate-based PFR reserve limits& Equivalency Ratio and Rate-Based PFR reserve Limits\\
(Optimization A or Opt. A)&(Optimization B or Opt. B)& (Optimization C or or Opt. C)\\
\hline

$\begin{aligned}[t]
\underset{b\in\mathbb{R}^{\beta}}{\underset{G\in\mathbb{R}^n,R\in\mathbb{R}^n}{\text{min}}} \  &\underset{i\in \mathbb{N}}{\Sigma}  C_{i}(G_{i}) \\[-3\jot]
& \ \text{st: } (\ref{PrivateBatteries}), (\ref{PrivateGenerators}), (\ref{PengweDuReq}), \text{ and }(\ref{LineConst})\nonumber
\end{aligned}$&$\begin{aligned}[t]
\underset{b\in\mathbb{R}^{\beta},r\in\mathbb{R}^n}{\underset{G\in\mathbb{R}^n,R\in\mathbb{R}^n}{\text{min}}} \  & \underset{i\in \mathbb{N}}{\Sigma}  C_{i}(G_{i})\\[-3\jot]
& \ \text{st: } (\ref{ResReqMain}), (\ref{PrivateBatteries}), (\ref{PrivateGenerators}), (\ref{ResReqRamp}), \text{ and }(\ref{LineConst})\nonumber
\end{aligned}$&$\begin{aligned}[t]
\underset{b\in\mathbb{R}^{\beta},r\in\mathbb{R}^n}{\underset{G\in\mathbb{R}^n,R\in\mathbb{R}^n}{\text{min}}} \  & \underset{i\in \mathbb{N}}{\Sigma}  C_{i}(G_{i})\\[-3\jot]
& \ \text{st: } (\ref{ResReqMain}), (\ref{PrivateBatteries}), (\ref{PrivateGenerators}), (\ref{ResReqRamp}), (\ref{EquivRatioPFRLimit}), \text{ and }(\ref{LineConst})\nonumber
\end{aligned}$\\
\hline	
\end{tabular}
\end{center}
\vspace{-20pt}
\end{table*}
\endgroup

\vspace{-2pt}
\subsection{Co-optimization Formulations}
\vspace{-4pt}
Without loss of generality, this section assumes that one generator is located at each bus in the system.  The total number of transmission lines is denoted $\ell$.  We use the DC approximation of the transmission network as outlined in~\cite{garcia2019approximating}.  The transmission constraints are as follows, where $S\!\in\!\mathbb{R}^{\ell\times n}$ is a DC approximation matrix, the line flow limits are denoted $\bar{F}\in\mathbb{R}^{\ell}$, and the fixed demand vector is denoted $D\in\mathbb{R}^n$:
\begin{equation}
-\bar{F} \leq S(G-D)\leq \bar{F} \ \ \ \text{and} \ \ \ \bold{1}^{\dagger}(G-D)=0\label{LineConst}
\end{equation} 
Table \ref{Table:OptProbs} provides the co-optimization problems, which differ only by their reserve requirements.  Opt. A enforces the equivalency ratio requirement (\ref{PengweDuReq}) from~\cite{liu2018participation} as detailed in Section~\ref{sec:EqRatReq} with parameters from Table \ref{DuTable}.  Opt. B and C enforce the general requirement (\ref{ResReqMain}) and the rate-based PFR reserve limit (\ref{ResReqRamp}) with governor ramp rates fixed to $\kappa_i\!=\!20$MW/s.  \mbox{Opt. C} also enforces the equivalency ratio PFR reserve limit (\ref{EquivRatioPFRLimit}). The proportional PFR reserve limit (\ref{Constraint:ConstTime}) is not directly studied in this section but is qualitatively similar to (\ref{EquivRatioPFRLimit}).

Opt. B and C are non-convex problems because the function $h(M,\bold{1}^{\dagger}b)$ is strictly convex in $b$ and appears on the Right-Hand-Side (RHS) of constraints (\ref{ResReqRamp}).  To formulate a convex problem, we suggest approximating the rate-based PFR reserve limit (\ref{ResReqRamp}) by evaluating the limit at an approximated value of total FFR reserve $\hat{b}\!\in\!\mathbb{R}^{\beta}$.  This allows the PFR reserve limit constraint to become linear in the optimization variables, e.g. $r_i\leq \kappa_i h(M,\hat{b})$ where the RHS is now constant.  As a result, Opt. B and C become linearly constrained and convex. Intuitively, this approximation uses estimates of total FFR reserve to systematically determine fixed PFR reserve limits.   

In the special case where FFR reserve is offered as zero price, Opt. B and C reduce to linearly constrained problems.  This is because it is always optimal to allocate FFR reserve at its maximum value, e.g. $b\!=\!\bar{b}$.  To solve Opt. B and C, we fix \mbox{$b\!=\!\bar{b}$} and solve the resulting linearly constrained convex problems.

\vspace{-15pt}
\subsection{Numerical Results}
\begin{figure}[h]
\centering
\vspace{-5pt}
\includegraphics[scale=.25]{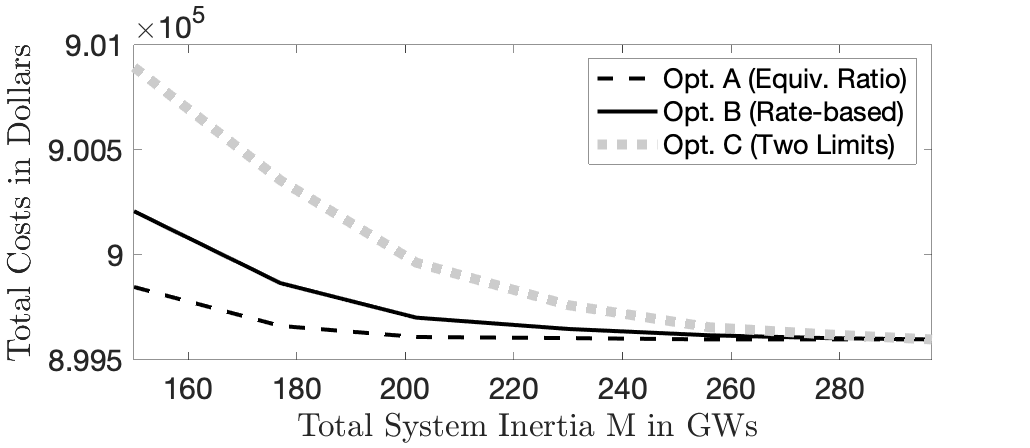}
\vspace{-15pt}
\caption[]{\label{fig:CoOptCosts} Optimal costs for each co-optimization formulation.}
\vspace{-15pt}
\end{figure}

\vspace{-3pt}
For each co-optimization problem, Figure \ref{fig:CoOptCosts} plots the total system cost as the fixed total system inertia $M$ decreases.  Each point plotted corresponds to an inertia level from Table~\ref{DuTable}. Notice that all formulations have the same objective value for $M\!=\!297$GWs.  This is because none of the reserve constraints are binding at this high inertia level. As the inertia decreases, the lowest costs result from Opt. A.  This is because the equivalency ratio requirement does not place a direct limit on the amount of PFR reserve that can be procured from a single generator and so significant amounts of low opportunity cost PFR reserve are procured from relatively few generators.  In contrast, \mbox{Opt. B} results in higher costs because it limits the procurement of low opportunity cost PFR reserve, forcing more expensive PFR reserve to be procured.  The cost for Opt. C is even higher than Opt. B because it enforces an additional constraint (\ref{EquivRatioPFRLimit}), reducing the size of its feasible set.  

\begin{figure*}[t!]
\vspace{-21pt}
\subfloat[Optimization A (Equivalency ratio requirement)]{\label{Fig:OptA}\includegraphics[scale=.21]{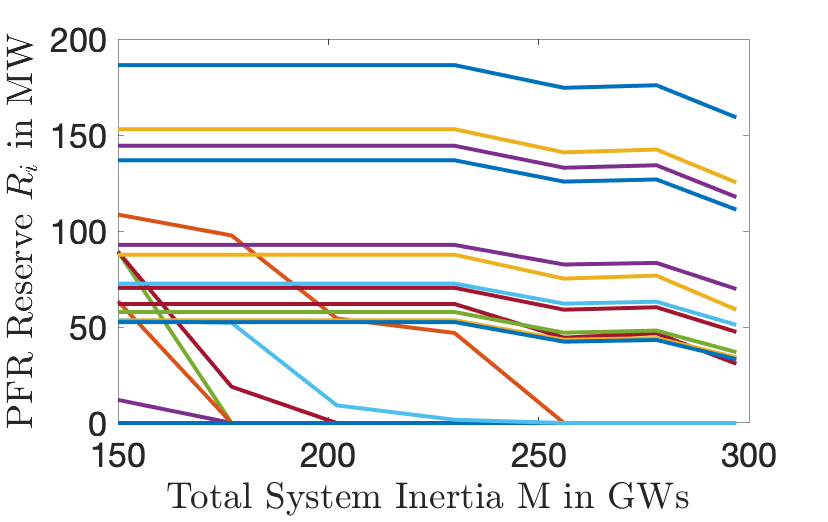}}
\subfloat[Optimization B (Rate-based PFR reserve limits)]{\label{Fig:OptB} \includegraphics[scale=.21]{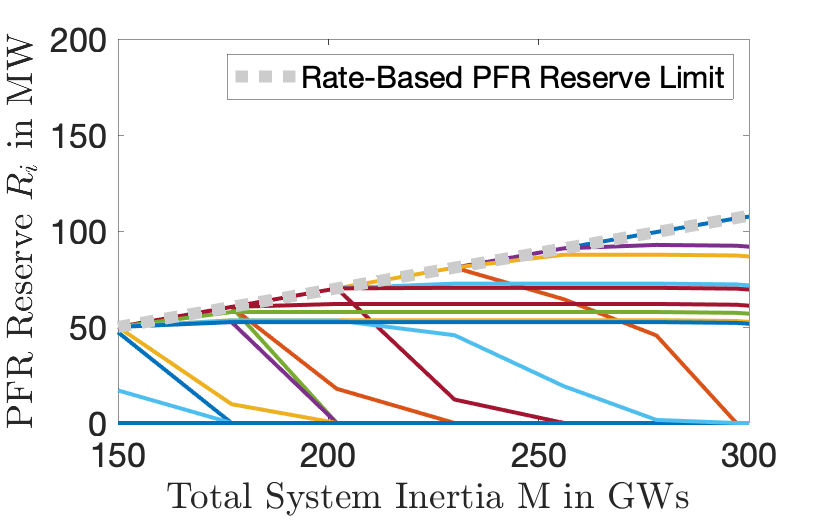}}
\subfloat[Optimization C (Two PFR reserve limits)]{\label{Fig:OptC} \includegraphics[scale=.21]{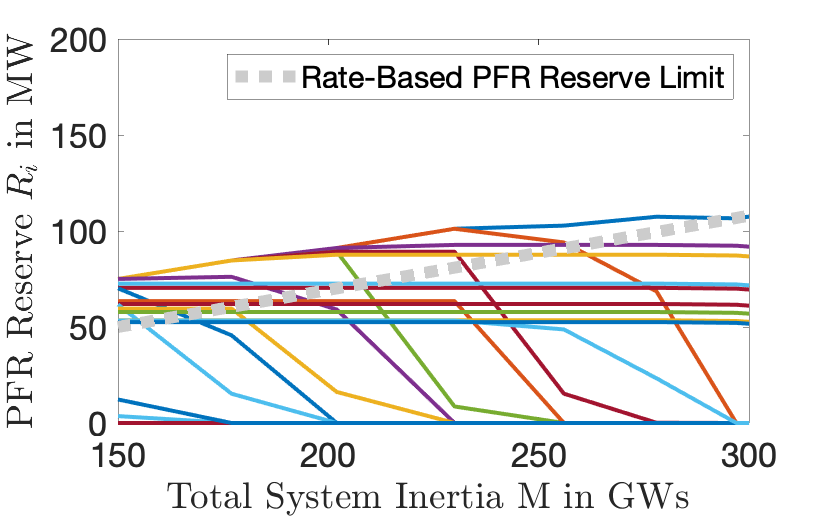}}
\caption{\label{fig:PFRProcurment}Three plots represent Optimization A, B, and C.  Each plot illustrates the procured nominal PFR reserve for each generator as inertia declines.  Each trajectory represents a different generator providing PFR reserve.}
\vspace{-18pt}
\end{figure*}

Figure \ref{fig:PFRProcurment} provides a plot for each co-optimization problem formulation that shows the procured nominal PFR reserve for each generator as the inertia $M$ decreases.  In this context, each trajectory represents the procured nominal PFR reserve for an individual generator. Each formulation has multiple solutions when the inertia is $297$GWs because the reserve constraints are not binding. As the inertia decreases, the reserve constraints become binding for each formulation.  The reserve constraints become binding for Opt. A, B, and C at $M\!=\!202$GWs, \mbox{$M\!=\!256$GWs,} and $M\!=\!278$GWs respectively.  

Let's first analyze Opt. B in Figure \ref{Fig:OptB}, which additionally shows the rate-based PFR reserve limit \mbox{$\kappa_i h(M,\bold{1}^{\dagger}\bar{b})$} as a dashed gray trajectory.  As the inertia drops, the limit function $h(\cdot,\bold{1}^{\dagger}\bar{b})$ decreases and the rate-based PFR reserve limit becomes binding for some generators.  This limit constraint causes the procured nominal PFR reserve to decrease for many generators with low procurement costs and increase for many generators with high procurement costs.  When the inertia reaches its lowest value of $152$GWs, many generators' procured nominal PFR reserve matches the rate-based PFR reserve limit.  It's apparent that Opt. B accommodates low inertia by distributing PFR reserve more evenly among generators.    

Figure~\ref{Fig:OptA} illustrates that Opt. A accommodates low inertia values by simply increasing the procured nominal PFR reserve for the generators with low reserve procurement costs.  Notice that all nominal PFR reserve trajectories increase as the inertia decreases in Figure~\ref{Fig:OptA}, which contrasts with Figure~\ref{Fig:OptB}.  The total procured reserve (including FFR and PFR) is equal to $L=2500$MW at all inertia values in Figure~\ref{Fig:OptB}, whereas the total procured reserve is more than $L=2500$MW in Figure~\ref{Fig:OptA}.

Figure~\ref{Fig:OptC} represents Opt. C, which intends to exhibit features of both the equivalency ratio requirement and the rate-based PFR reserve limit. Similar to Opt. B, Opt. C sees many generators providing nominal PFR reserve at low inertia levels and sees no single generator offering a substantial amount of PFR reserve.  Similar to Opt. A, Opt. C leaves extra headroom for each generator, thus the nominal PFR reserve for some generators is greater than the rate-based PFR reserve limit and the total nominal PFR reserve is more than $L=2500$MW.

\vspace{-12pt}
\section{Conclusions} \label{Sec:Conc}
\vspace{-5pt}

This paper derives reserve requirements that couple PFR reserve, FFR reserve, and inertia.  The proposed requirements state that sufficient FFR and PFR reserve must be procured to cover the largest possible loss of generation and that this reserve is available to be deployed before the frequency reaches the critical frequency threshold.  The amount of available PFR reserve provided by a generator is limited by its ramping ability.  We derive two such PFR reserve limits from first principles termed the \emph{rate-based PFR reserve limit} and the \emph{proportional PFR reserve limit}.  We additionally derive the \emph{equivalency ratio PFR reserve limit} empirically from the \emph{equivalency ratio reserve requirement}.  These limits are used, along with various assumptions, to derive the equivalency ratio from first principles, which has only been studied empirically to date.  Numerical results illustrate the high-level differences between each reserve requirement as total system inertia decreases.  The equivalency ratio requirement concentrates large amounts of PFR reserve to relatively few generators at low inertia levels.  In contrast, the proposed requirement disperses the procured PFR reserve more evenly among many generators when enforcing the rate-based PFR reserve limit.  When enforcing the equivalency ratio PFR reserve limit, our requirement captures qualities of the equivalency ratio requirement.

\vspace{-15pt}

\bibliography{bibfile}
\bibliographystyle{IEEEtran}

\end{document}